\documentclass[preprint,12pt]{elsarticle}
 \usepackage[T1]{fontenc}
\usepackage{amsmath}
\usepackage{graphicx}
\usepackage{mathtools}   % loads »amsmath«
\usepackage{amssymb}

\usepackage{stmaryrd}
\usepackage{bm}
\usepackage{amsfonts} 
\usepackage{float}
\usepackage{xcolor}
\usepackage{amsmath,bm}
\usepackage{mathabx}
\newcommand*\xbar[1]{%
  \hbox{%
    \vbox{%
      \hrule height 0.5pt % The actual bar
      \kern0.5ex%         % Distance between bar and symbol
      \hbox{%
        \kern-0.1em%      % Shortening on the left side
        \ensuremath{#1}%
        \kern-0.1em%      
      }%
    }%
  }%
} 
\newtheorem{theorem}{Theorem}[section]

\newtheorem{proposition}[theorem]{Proposition}

\newtheorem{definition}[theorem]{Definition}%{\begin{trivlist}

\usepackage{amssymb}
\usepackage[nodots]{numcompress}
\makeatletter
\newcommand{\ostar}{\mathbin{\mathpalette\make@circled\star}}
\newcommand{\make@circled}[2]{%
  \ooalign{$\m@th#1\smallbigcirc{#1}$\cr\hidewidth$\m@th#1#2$\hidewidth\cr}%
}
\newcommand{\smallbigcirc}[1]{%
  \vcenter{\hbox{\scalebox{0.77778}{$\m@th#1\bigcirc$}}}%
}
\begin{document}

\begin{frontmatter}

\tnotetext[label1]{Corresponding author: Terence R. Smith, Department of Computer Science, University of California at Santa Barbara, Ca 93106}
\author{Terence R. Smith { }\fnref{label2}}
 \ead{smithtr@cs.ucsb.edu}
\title{\bf Sums of Exponential Terms, Conserved Quantities,
and the Real Wave Numbers}
\begin{abstract}
There is consensus that sums 
$S_n={ {\Sigma }_{k=1}^n R_{0k} e^{i \theta_k}}$
of complex exponential terms,
despite their mathematical significance,
only possess closed-form representations 
for specific values of n and special values of their parameters
and that there are no generally-accepted recursive formulae for 
their computation.
This note is focused on
recursive formulae that: 
(1) provide closed-form analytic representations of $S_n$
for any finite n; 
(2) include generalizations of the usual formula
for the sum of two exponentials; and 
(3) are representable in the form
$S_n= A_n exp({  i\Sigma_{k=1}^n \theta_k})$.
The goal of the paper is to show that one may interpret 
the exponential term
$exp(i \Sigma_{k=1}^n \theta_k)$
of $S_n$ as representing the projection,
from a field of numbers that generalizes the complex numbers
onto the complex plane, of a term representing quantities that are conserved
under the addition and multiplication
of numbers in the extended space.
In particular, it is shown that the general form of a number in the extended field
generalizes the form of a sum of complex exponentials.
\end{abstract}

\begin{keyword}
\ Sums of complex exponentials, canonical representation, real wave numbers, conserved quantities
\end{keyword}
 
\end{frontmatter}
%\vspace{10px}

%%%%%%%%%%%%%%%%%%%%%%%%%%%%%%%%%%%%%%%%%%%%%%%%%%%%%%%%

\section{Introduction}

Linear combinations of exponential terms 
\begin{eqnarray}
\label{eq:series}
S_n={\Sigma }_{k=1}^n R_{0k} e^{i \theta_k}, \hspace{.5in} \{R_{0k}, \theta_k | k=1,n\}\  \epsilon\ \mathbb{C}
\end{eqnarray}
are important, with applications that include
the representation of periodic functions [1] and representations of solutions to linear differential equations [2].
A search of the mathematical literature suggests a consensus 
that: 
(1) the general sum $S_n$ does not possess a closed-form representation;
(2) closed-form representations may only be found for
specific values of n
or for cases in which the arguments $\{R_{0k}, \theta_k | k=1,n\}$
take special forms [3];
and (3) there is no
generally-accepted recursive formula for computing the general sums $S_n$,
other than appling Euler's relation $e^{i\theta}=\cos{\theta}+i \sin{\theta}$
to partial sums $S_{n+1}= S_n + R_{0 {n+1}}e^{i \theta_{n+1}}$
of exponential terms.

It is shown in this note that there are
recursive formulae that not only provide closed-form analytic representations
of $S_n$ for any n, but also possess a canonical form that
involves the exponential term $exp({  i\Sigma_{k=1}^n \theta_k})$ of $S_n$.
It is also shown that the canonical nature of this form is made explicit 
in the arithmetic of a field of numbers that extends the
complex numbers $\mathbb{C}$ to exponential functions of linear mappings
of the real numbers. In particular, it is shown that the sum $S_n$
has a form that is analogous to 
the form of the numbers in this space, 
and that the argument $ \Sigma_{j=1}^n \theta_j$
of its exponential term
is a projection onto the complex plane of one of two quantities
that are conserved in the arithmetic of the extended space of numbers.

Two facts facilitate the analysis. 
First,
sums of two exponentials may be written as
\begin{eqnarray}
\label{lemma}
e^{i\theta_1} +e^{i\theta_2}=
\Big(e^{i\big(\tfrac{\theta_1-\theta_2}{2}\big)}+e^{-i\big(\tfrac{\theta_1-\theta_2}{2} \big)}\Big)
e^{i\big(\tfrac{\theta_1+\theta_2}{2}\big)}
=2cos\big(\tfrac{\theta_1-\theta_2}{2}\big)e^{i\big(\tfrac{\theta_1+\theta_2}{2}\big)} 
\end{eqnarray} 
\noindent 
in which the first equality is an identity 
and the second is a definition of the cosine.
Second, one may assume that coefficients $\{R_{0k}, k=1,n \}$
of Equation (\ref{eq:series})
are non-negative real numbers. Otherwise,
any complex-valued coefficient $R_{0k}$ may be multiplied
by $|{R}_{0k}|/|{R}_{0k}|$, in which $|\bar{R}_{0k}|$
is the magnitude of $R_{0k}$,
and $R_{0k}/|{R}_{0k}|$ 
represented as
$e^{i\phi_k}$ for some $\phi_k$,
leaving the form  of the sum (\ref{eq:series}) unchanged.

\section{Recursive Formulae for Sums of Exponential Terms}

There are recursive representations of the sum $S_n$ in which
the term $exp(i {\Sigma }_{k=1}^n\theta_k)$ occurs.
An informative representation follows from
an identity of $S_n$
for $\{R_{0k}=1|k=1,n\}$:
\begin{eqnarray}
\label{eq:A0}
{\overset{n} {\underset {j=1} \Sigma }}e^{i \theta_j}
=\Big(\overset {n} {\underset {j=1} \Sigma } 
e^{-i\overset {n} {\underset {k\ne j} \Sigma }(\theta_k)}\Big)\ 
e^{i\overset {n} {\underset {j=1} \Sigma }(\theta_j)}
\equiv A_ne^{i\overset {n} {\underset {j=1} \Sigma }\theta_j}
%= A_n e^{i\overset {n} {\underset {j=1} \Sigma }\theta_j}
\end{eqnarray}
which may be loosely interpreted as stating that the term $A_n$
of Equation (\ref{eq:A0})
represents the factoring of non-conjugate exponential terms
from the sum on the LHS of  (\ref{eq:A0}).
When generalized with amplitudes $R_{0k}$, Equation (\ref{eq:A0}) becomes
\begin{eqnarray}
\label{eq:A1}
{\overset{n} {\underset {j=1} \Sigma }}R_{0j} e^{i \theta_j}
={\overset{n} {\underset {j=1} \Sigma }} e^{i (-ilnR_{0j}+\theta_j})
=\Big(\overset {n} {\underset {j=1} \Sigma } 
e^{-i\overset {n} {\underset {k\ne j} \Sigma }(-i lnR_{0k} +\theta_k)}\Big)\ 
e^{i\overset {n} {\underset {j=1} \Sigma }(-i lnR_{0j} +\theta_j)}
\end{eqnarray}
which representats the sum of n terms
as n sums of $n-1$ terms. In summary, one has
\begin{proposition}
\label{prop:1}
The sum
of n exponential terms
with $R_{0k} > 0$ 
may be represented as
\begin{eqnarray}
\label{eq:A3}
{\overset{n} {\underset {j=1} \Sigma }} 
R_{0j} e^{i \theta_j}
= A_n e^{i\overset {n} {\underset {j=1} \Sigma }\theta_j},\ \  \hspace{.2in} 
 A_n=\big(\overset {n} {\underset {j=1} \Pi }R_{0j}\big)
\Big(\overset {n} {\underset {j=1} \Sigma } 
e^{-i\overset {n} {\underset {k\ne j} \Sigma }(-i ln(R_{0k} +\theta_k)}\Big) 
\end{eqnarray}.
\end{proposition}
\noindent
There are various alternative representations of the term 
$A_n$. A first is stated in
\begin{proposition}
\label{thm:wnform11}
\hspace{.2in} $A_{n}=
A_{n-1}e^{-i \theta_{n} }+R_{0n}e^{-i\overset {n-1} {\underset {j=1} \Sigma }\theta_j}, \hspace{.2in} A_0=0$\\
\end{proposition}
\noindent
whose proof follows upon substituting the recursive relation for  $A_n$
into Equation (\ref{eq:A3}) to obtain 
the recursive identity
 $S_{n+1}= S_n + R_{0 {n+1}}e^{i \theta_{n+1}}$.

A third representation for $A_n$
that generalizes the form given for a sum of two exponential terms
in Equation (\ref{lemma})
may be found by applying
the second equality of Equation (\ref{lemma}) to the definition
of $A_n$ in Proposition 2 and simplifying, which leads to 
\begin{proposition}
\label{prop2}
\hspace{.35in}$A_{n}
=(A_{n-1}R_{0n})^{1/2}
cos\Big(   \tfrac{  
-iln\big(\frac{A_{n-1} } {R_{0n}}\big)+
 {{\overset {n} {\underset {j=1} \Sigma }\theta_j}}-\theta_{n}}{2}\Big)
e^{-i\tfrac{\overset {n} {\underset {j=1} \Sigma }\theta_j}{2} }$
\end{proposition}
\noindent
One notes that this representation of $A_n$ involves 
an exponential term and that the recursive structure 
of the argument of the cosine function
implies that the cosine terms
are nested, with an additional level of nesting for every additional exponential term in the sum $S_n$, so giving rise to $n$ levels of nesting in the cosine terms.

One may construct representations of $S_n$
that do not include the term $exp(i {\Sigma }_{k=1}^n\theta_k)$:
\begin{proposition}
\label{thm:wnform12}
The sum
of n exponential terms
with $R_{0k}\geq 0$
may be represented as
\begin{eqnarray}
\label{eq:cor1}
S_n &=&\overset {n} {\underset {j=1} \Sigma }
R_{0j} e^{i \theta_j}
= A_n e^{i\sigma_n}\hspace{.1in} \ \ \\
in\
 which\ \hspace{.1in} \sigma_n&=&\tfrac{\theta_1 + \overset {n} {\underset {j=2} \Sigma }2^{j-2} \theta_j}
{2^{n-1}},\ \ \ \ \ \ \sigma_1=\theta_1; \\
A_n&=&\begin{cases}
\label{eqn}
    R_{01}, \hspace{3in} n=1\\
    (A_{n-1}R_{0n})^{1/2} 2 cos(-i ln \big(\tfrac{A_{n-1}}{R_{0n}}\big)^{1/2} + \tfrac{1}{2} (\sigma_{n-1} - \theta_n)), \ n\ge 2
  \end{cases}
\end{eqnarray}
\end{proposition}
\noindent
which has a straightforward inductive proof.
This representation may, however, be transformed in an informative manner
into the general form
represented in Propositions 1-3. 
By the commutativity 
and associativity of complex arithmetic, the representations of the sum $S_n$
in Proposition 1-4
are symmetric in value
under permutations of the arguments ($R_{0k}, {\theta}_k)$ for $ k=1,n$,
hence
on taking the nth root
of the product of all cyclic permutations
$\{P_n^{k-1}, k=0, n-1\}$ of the sum $S_n$, one obtains
\begin{eqnarray}
S_n 
&=& 
\big(\overset {n} {\underset {j=1} \Pi }   P_n^{k-1} \big(S_n\big)\big)^{1/n}
=
\big(\overset {n} {\underset {j=1} \Pi }   P_n^{k-1} \big(A_n\big)\big)^{1/n}
\big(\overset {n} {\underset {j=1} \Pi }   P_n^{k-1} \big(e^{i\theta_n}\big)\big)^{1/n}
\end{eqnarray}
and notes that the exponential term 
$exp\big({i \Sigma_{j=1}^n\theta_j}\big)$ of Propositions 1-3 are invariant in both form and value. 
While the corresponding term
of Propositions 4 is not invariant and takes the form 
$exp\big({i \Sigma_{j=1}^n\theta_j/n}\big)$,
one notes that $exp(i {\Sigma }_{k=1}^n\theta_k/n)
=exp(-i {\Sigma }_{k=1}^n\theta_k (1-1/n).exp(i {\Sigma }_{k=1}^n\theta_k)$.
Hence the invariant form may be restored by absorbing the term $ exp(-i {\Sigma }_{k=1}^n\theta_k (1-1/n)$ into the coefficient $A_n$.

A natural question that arises from these observations concerns
the significance of the exponential term 
$exp(i {\Sigma }_{k=1}^n\theta_k)$.
In particular, the invariance of its form and value under permutations
suggests the existence of an invariant, or conserved, quantity.

\section{A GENERALIZATION OF THE COMPLEX NUMBERS}
In seeking to answer the preceding question, it is of value
to consider 
a space of elements that are generated from the set of exponential mappings 
of the elements of the vector space L($\mathbb{R},\mathbb{R}$)
of linear mappings from $\mathbb{R}$ onto 
$\mathbb{R}$. 
Generators for this space of elements
may be defined in terms of the set of functions
\begin{eqnarray}
\label{eq:multi}
{\bf w}(f, \theta)=\{e^{i 2 \pi (f{\bm \rho}+\theta)}| \ 
\forall \rho\ \epsilon\ \mathbb{R}\},
\forall\ f,\theta 
\ \epsilon\  \mathbb{R}.
\end{eqnarray}
with each function having a period $1/f$,
a translation $\theta$,
and an interpretation as an infinite-dimensional extension of the complex number
$e^{i2\pi \theta}$.
A discussion of wave numbers for the case in which $\rho$
assumes integer values is presented in [5].

It is natural to ask whether one may 
multiply, divide, add, and subtract the elements 
${\bf w}(f, \theta)$ and whether they form a field of numbers
under the closure of these operations.
It is straightforward to show that this is the case if one defines 
operators 
for multiplication $\otimes$, addition $\oplus$, and inverse $I$
to be operators with the usual definitions of $\{+,-,x,\div\}$
applied in a pointwise manner.
One notes that the pointwise 
application of these operators ensures that the usual associative,
commutative, and distributive laws of arithmetic hold
when applying the operators $\otimes, \oplus$, and $I$
to the elements.

The generators of Equation (\ref{eq:multi}) form an Abelian
multiplicative group under $\otimes$ with
\begin{eqnarray}
{\bf w}(f_1, \theta_1) \otimes {\bf w}(f_2, \theta_2)=
{\bf w}(f_1+f_2, \theta_1+\theta_2)
\end{eqnarray}
since, by the existence of the inverse $I$, one has $I({\bf w}(f, \theta))
=\overline{{\bf w}}(f, \theta)$, in which $\overline{\bf w}$
denotes complex conjugation, and hence the multiplicative identity
${\bf w}(0, 0)$.
In showing that they generate an Abelian additive group under $\oplus$,
one may apply any of Propositions 1-3 in a pointwise manner
to the elements of the multiplicative group
to obtain
\begin{eqnarray}
{\bf w}(f_1, \theta_1) \oplus{\bf w}(f_2, \theta_2)
&=&
 {\bf A}_{12}{\bf w} ({ f_1+f_2,  \theta_1+\theta_2})\\
\label{times}
{\bf A}_1{\bf w}(f_1, \theta_1) \oplus {\bf A}_2{\bf w}(f_2, \theta_2)
&=&
e^{ i (-iln (   {\bf A}_1) \oplus     
(f_1{\bm \rho}+\theta_1 ) }     \oplus          
e^{ i (-iln (  {\bf A}_2) \oplus     
(f_2{\bm \rho}+\theta_2)           }\\
&=&
 \hat{{\bf A}}_{12}{\bf w} ({ f_1+f_2,  \theta_1+\theta_2}) \nonumber
\end{eqnarray}
in which ${\bf A}_{1},{\bf A}_{2}, {\bf A}_{12},  \hat{{\bf A}}_{12}$ may be written 
as sequences of trigonometric
and exponential terms.
Since the set of elements is closed under $\otimes$, 
multiplication of an element ${\bf A} {\bf w}(f, \theta)$ by ${\bf w}(0, 1/2)$ leads to 
$-{\bf A} {\bf w}(f, \theta)$, which is its additive inverse,
and hence to the additive identity 
${\bf 0}={\bf A} {\bf w}(f, \theta)\oplus(-{\bf A} {\bf w}(f, \theta))$. 
Furthermore, the existence of the inverse of an element
$I({\bf A} {\bf w}(f, \theta))
=\big(1/{\bf A} \big) \overline{{\bm \omega}}(f, \theta)$ 
leads to the multiplicative identity and to the fact
that the elements form a multiplicative group.
Since pointwise division $\odiv$ is defined by the product of one element
with the inverse of another, 
it leads to a field of elements
that may be termed the real wave numbers
and denoted by $\mathbb{W}$.

It is clear that the application of any of the operators $\otimes, \oplus, I$
leaves the form
\begin{eqnarray}
\label{eq:A2}
{\bm \omega} = {\bf A} {\bf w} (f,\theta)
\end{eqnarray}
invariant, which remains the case under the closure of the operators.
In particular, since $\bar{\bf w}(f,\theta)\epsilon \mathbb{W}$, all products,
sums, and inverses of the terms ${\bf A}$ are elements of $\mathbb{W}$
and these may be represented as ${\bm \omega} = {\bf A} {\bf w} (0,0)$.
From the general form of a product
\begin{eqnarray}
\label{prod}
{\bf A}_1{\bf w}(f_1, \theta_1) \otimes {\bf A}_2{\bf w}(f_2, \theta_2)
&=&
{\bf A}_{1}{\bf A}_{1}{\bf w} ({ f_1+f_2,  \theta_1+\theta_2})
\end{eqnarray}
and the analogous form for a sum in Equation (\ref{times}), 
one notes that dilations
and translations are conserved quantities under 
summation and multiplication.

It is of interest to ask whether there is a proper subfield 
of $\mathbb{W}$ of periodic functions,
and this too may be answered in the affirmative. While a product of numbers 
of the form ${\bf w}(f,\theta)$ is always periodic, it is generally the case
that a finite sum of such numbers
is periodic if and only if the ratios of 
their periods $\{f_j/f_k | j,k =1,n\}$ are rational numbers [4].
It follows that while wave numbers ${\bm \omega} \epsilon {\mathbb{W}}$
are generally not periodic,
they are periodic in the proper subfield in which
$f,\theta$ are rational numbers.

\section{INTERPRETING THE FORM OF SUMS OF EXPONENTIAL TERMS}
One notes that the general form (\ref{eq:A2}) of a real wave number
is analogous to
the form $ A_n e^{  i\Sigma_{k=1}^n \theta_j}$ of a sum
of exponentials, as defined in Propositions 1-3,
in terms of its
generalized amplitude ${\bf A}_n$
and generalized exponential term ${\bf w}(f,\theta)$. This
suggests that wave numbers may provide useful insights
into sums of exponential terms.

One such insight follows on noting that increasing
(or decreasing) values of the wave number parameter ${\bm \rho}$
may be viewed as moving a point representing 
a wave number across sheets defined above and below
the complex plane, as for example, in the case of the complex logarithmic
function. Such motion may be viewed as occurring in either 
a clockwise or counter clockwise direction relative to the unit circle,
and a wave number viewed as tracing out the form
of a compressed, generalized helix.

Given this interpretation of the real wave numbers,
it is natural to make
\begin{definition}
The spin and rotation of a real wave number
${\bf A} {\bf w}(f,\theta)$
are $f$ and $\theta$.
\end{definition}
\noindent
This definition, together with Equations (\ref{times}) and (\ref{prod}),
then leads to
\begin{proposition}
Both the spin and rotation of real wave numbers are
additively and multiplicatively conserved quantities.
\end{proposition}
\noindent
The form of the exponential term ${\bf w}(f,\theta)$ of the real wave numbers 
and its associated 
additive and multiplicative invariants 
together suggest
that the corresponding exponential
term $exp(i {\Sigma }_{k=1}^n\theta_k)$
occurring in representations of sums of ordinary exponentials 
may be interpreted as a projection onto the complex plane
of the spins and rotations of the wave numbers
in their various sheets.
This interpretation is also suggested by the nested cosine terms
of Equation \ref{prop2}. 

The form of the exponential term in Propositions 1-3
may therefore be viewed as canonical
in the sense that it represents the projection
of conserved quantities from a more general space of numbers.
The projection is not lossless, since information about the spin
is lost, while only the information concerning the rotation is preserved.

This interpretation of the exponential term
in the expression for the sum $S_n$
suggests that other properties of numbers
in the complex plane may be more easily understood
in terms of the properties of real wave numbers.
This is analogous to the observation that
there are properties of numbers
on the real line that are sometimes more easily understood
in terms of projections from the complex plane
[6].
\section{Bibliography}
\noindent
[1]
Evans, L.C., 2002.
{Partial Differential Equations},
2nd  edition, American Mathematical Society, Washington DC. 
{662} pages.\\
\noindent
[2] Goodman, J W, {2015}.
{Statistical Optics},
{2nd} edition,
 {Wiley},  {Hoboken New Jersy},
{544} pages.\\
\noindent
[3]
{ Stein, E M and R Shakarchi}, 2003.
{Statistical Fourier Analysis: An Introduction},
1st  edition,
 {Princeton University Press},
{Princeton New Jersy},
328 pages.\\
\noindent
[4]
{Mirotin, A R and E A Mirotin}, 2009.
 {On Sums and Products of Periodic Functions},
 {Real Analysis Exchange},
  volume	 {34},
  number {2},
  pages{347-358}.\\
\noindent
[5]
{Smith, T R}, 2025.
 {Rational Wave Numbers and the Algebraic 
Structure of the Cyclic Groups of the
Roots of Unity},
 { arXiv:2503.07629v1 [math.NT]}.\\
\noindent
[6]{ Needham, T}, 2023.
{Visual Complex Analysis},
{Anniversary   edition},
{Oxford University Press},
 {Oxford, England},
{328} pages.

%\bibliographystyle{vancouver}
%\bibliography{VancouverExamples.bib}
\end{document}